\documentclass[11pt,twoside]{article}
\usepackage{amsmath,amsfonts,amssymb,theorem}

\def\lhead{K.Kang }
\def\rhead{Partial regularity of the Navier-Stokes equations}
\def\Xint#1{\mathchoice
   {\XXint\displaystyle\textstyle{#1}}%
   {\XXint\textstyle\scriptstyle{#1}}%
   {\XXint\scriptstyle\scriptscriptstyle{#1}}%
   {\XXint\scriptscriptstyle\scriptscriptstyle{#1}}%
   \!\int}
\def\XXint#1#2#3{{\setbox0=\hbox{$#1{#2#3}{\int}$}
     \vcenter{\hbox{$#2#3$}}\kern-.5\wd0}}

\def\aint{\Xint\diagup}

\newcommand{\e}{\varepsilon}

\newcommand{\raroup}{\rightharpoonup}

\newcommand{\bke}[1]{\left( #1 \right)}
\newcommand{\bkt}[1]{\left[ #1 \right]}
\newcommand{\bket}[1]{\left\{ #1 \right\}}
\newcommand{\norm}[1]{\left\Vert #1 \right\Vert}
\newcommand{\abs}[1]{\left| #1 \right|}

\newcommand{\bkp}[1]{\big( #1 \big)}

\newcommand{\R}{{\mathbb R }}
\newcommand{\calS}{{\cal S }}
\newcommand{\ap}{{\alpha}}
\newcommand{\calC}{{\mathcal C}}
\newcommand{\calO}{{\mathcal O}}
\newcommand{\calP}{{\mathcal P}}
\newcommand{\pr}{{\partial}}
\newcommand{\om}{{\Omega}}
\newcommand{\eps}{{\epsilon}}
\newcommand{\nb}{{\nabla }}
\newcommand{\ind}{{\,\mathrm{d}}}
\newcommand{\raro}{{\rightarrow}}
\newcommand{\be}{\begin{equation}}
\newcommand{\ee}{\end{equation}}
\newcommand{\bes}{\begin{equation*}}
\newcommand{\ees}{\end{equation*}}
\newcommand{\bea}{\begin{eqnarray}}
\newcommand{\eea}{\end{eqnarray}}
\newcommand{\ps}{{Prodi-Serrin }}

\newcommand{\sidenote}[1]{}

\newtheorem{thm}{Theorem}[section] {
  \newtheorem{assume}[thm]{Assumption}
  \newtheorem{lemm}[thm]{Lemma}
  
  \newtheorem{cor}[thm]{Corollary}
  \newtheorem{defi}[thm]{Definition}
{\theorembodyfont{\rmfamily}  \newtheorem{rem}[thm]{Remark} }
  
 }

\newenvironment{qed}{\hfill\fbox{}\par\vspace{.2cm}}

\newenvironment{pf}{{\bf Proof.}}{\hfill\fbox{}\par\vspace{.2cm}}
\newenvironment{lpf}{{\par\noindent\bf
            The sketch of the proof of Lemma \ref{mod:lemma} }}{}
\newenvironment{ispf}{{\par\noindent\bf
            The sketch of the proof of Lemma \ref{int:lemma100} }}{}
\newenvironment{sol}{{\par\noindent\bf
            The Proof of Theorem \ref{TH1} }}{}
\newenvironment{sol1}{{\par\noindent\bf
            The Proof of Theorem \ref{TH2} }}{}

\title{Regularity criteria for suitable weak solutions
of the Navier-Stokes equations near the boundary}

\author{ Stephen Gustafson, Kyungkeun Kang, and Tai-Peng Tsai }
\date{March 31, 2005}



\begin{document}
\maketitle

\begin{abstract}
We present some new regularity criteria for ``suitable weak
solutions'' of the Navier-Stokes equations near the boundary in
dimension three.  We prove that suitable weak solutions are H\"older
continuous up to the boundary provided that the scaled mixed norm
$L^{p,q}_{x,t}$ with $3/p+2/q\leq 2, 2<q\le \infty$, $(p,q) \not =
(3/2,\infty)$, is small near the boundary. Our methods yield new
results in the interior case as well. Partial regularity of weak
solutions is also analyzed under some conditions of the \ps type.
\end{abstract}

\section{Introduction}
In this note, we study the boundary regularity problem for {\it
suitable weak solutions} $(u,p):Q_T \to \R^3 \times \R$ to the
Navier-Stokes equations
\begin{equation}
\label{intro:nse}
\left\{
\begin{array}{c}
u_t-\Delta u+(u\cdot\nb)u+\nb p=f\\
\rm{div}\,\, u=0
\end{array}
\right.
\quad \mbox{ in }\,\,Q_T=\om\times (0,T),
\end{equation}
where $\om$ is a domain in $\R^3$, $u(x,t)$ is the velocity field and
$p(x,t)$ is the pressure.  By suitable weak solutions we mean
functions which solve the Navier-Stokes equations in the sense of
distribution, satisfy some integrability conditions, and satisfy the
local energy inequality (for details, see Definition
\ref{defsws:100} in section 2).  We assume that $\om$ and $f$ are
sufficiently regular and will give the specifics later.
For a point $z=(x,t)$ in $\overline \om \times (0,T]$,
denote $B_{x,r}= \{ y\in \R^3: |y-x|<r \}$,
\begin{equation*}
Q_{z,r} :=B_{x,r} \times (t-r^2, t) , \quad
Q^+_{z,r}:=Q_{z,r}\cap Q_T.
\end{equation*}
A solution $u$ is said to be {\it regular} at $z$ if $u$ is
uniformly H\"older continuous (for some exponent) in both $x$ and $t$
in $Q^+_{z,r}$ for some $r>0$.

After the seminal work of Leray \cite{Leray} and Hopf \cite{Hopf} on
the existence of weak solutions, the problems of uniqueness and
regularity of weak solutions remain unsolved.  A set of criteria which
guarantee uniqueness and regularity is the {\it \ps conditions}: any
weak solution $u$ of the Navier-Stokes equations is unique and regular
in $Q_T$ if it satisfies, for some $p,q \ge 1$,
\begin{equation} \label{eq:1-5}
 \norm{u}_{L^{p,q}(Q_T)} < \infty, \quad \frac 3p + \frac 2q \le 1 .
\end{equation}
Here
\begin{equation*}
   \norm{u}_{L^{p,q}(Q_T)} :=
   \norm{\norm{u(\cdot,t)}_{L^p_x(\om)}}_{L^q_t(0,T)}.
\end{equation*}
Note that a weak solution satisfying \eqref{eq:1-5} is automatically a
suitable weak solution: interpolating with $u \in L^{2,\infty}(Q) \cap
L^{6,2}(Q)$, a weak solution $u$ satisfying \eqref{eq:1-5} belongs to
$L^{4,4}(Q)$.  Hence one can use $u$ multiplied by a cut-off
function as a test function and derive the local energy inequality from
the weak formulation of \eqref{intro:nse}.  We now briefly review
regularity results for the \ps class. See \cite{Sohr01,ESS} for
more references.  Assuming \eqref{eq:1-5},
Serrin \cite{S62,S63} proved regularity when $3/p+2/q < 1$.  The cases
$3/p+2/q = 1$, $3<p\le \infty$ were proved by Fabes, Jones and
Rivi\'ere \cite{FJR} for $\om =\R^3$, by Sohr \cite{Sohr83} and Giga
\cite{G86} for $\om$ a domain, and by Struwe \cite{S88} for the
interior case.  See \cite{Sohr01,CP} for results in the setting of
Lorentz and Morrey spaces.  These results were recently extended up to
a flat boundary by the second author \cite{k_bnse} and to a curved
boundary by Solonnikov \cite{VS02}.  A flat boundary is a portion of
the boundary which lies on a plane. The endpoint case
$(p,q)=(3,\infty)$ was recently resolved by Escauriaza, Seregin and \v
Sver\' ak \cite{ESS} for the $\R^3$ and interior cases, and by Seregin
\cite{S04} for domains.

Recently, there have been many works on regularity criteria with
conditions involving only $p$. We will not try to give a list here.

After the partial regularity theory of Scheffer in a series of
papers \cite{VS76,VS77,VS80,VS82}, Caffarelli-Kohn-Nirenberg \cite{CKN}
proved that the one dimensional parabolic Hausdorff measure of the set
$\calS$ of possible interior singular points of suitable weak solutions
is zero, denoted $\calP^1(\calS)=0$.  This implies that the one
dimensional Hausdorff measure of $\calS$ is also zero. See section 2
for the definition of parabolic Hausdorff measures.
The key to the analysis in \cite{CKN} is the following regularity
criterion.  There is an absolute constant $\eps >0$ such that, if $u$
is a suitable weak solution of the Navier-Stokes equations in $Q_T$
and for an interior point $z=(x,t)\in Q_T$,
\begin{equation}
\limsup_{r\raro 0_+}\frac{1}{r}\int_{Q_{z,r}}\abs{\nb u(y,s)}^2dyds
\leq \eps,
\label{nse:abs10}
\end{equation}
then $u$ is regular at $z$. See \cite{L98} for a simplified
proof and \cite{LS99} for more details.

Recently, Seregin \cite{S02} extended the interior partial regularity
result up to a flat boundary.  More precisely, there exists an
absolute constant $\eps>0$ such that, if a suitable weak solution $u$
satisfies
\begin{equation}
\limsup_{r\raro 0_+}\frac{1}{r}\int_{Q^+_{z,r}}\abs{\nb u(y,s)}^2dyds
\leq \eps,
\label{nse:abs20}
\end{equation}
where $z\in \Gamma\times (0,T)$ and $\Gamma$ is a flat boundary of
$\om$, then $u$ is regular at $z$.
Combining the results in \cite{CKN} and \cite{S02}, one can conclude
that suitable weak solutions are H\"older continuous up to the
flat boundary away from a closed set $\calS\subset\overline Q_T$ with
$\calP^1(\calS)=0$.
The same assertion for a curved boundary is believed to be true,
but there seems no written  proof yet.

The objective of this paper is to present new sufficient conditions
for the regularity of suitable weak solutions to the Navier-Stokes
equations near the flat boundary (and in the interior).  Our main
result is that, in place of condition \eqref{nse:abs20}, H\"older
continuity of $u$ near the boundary can be ensured by the smallness of
the scaled mixed $L^{p,q}$-norm of the velocity field $u$.  We assume
that $f$ belongs to $M_{2,\gamma}$ for some $\gamma>0$ (this is 
a parabolic Morrey space, to be defined in section 2).  
We have the following theorem.

\begin{thm}[Regularity Criteria] \label{TH1}
Suppose $f\in M_{2,\gamma}(Q)$ for some $\gamma>0$, a parabolic Morrey
space.  For every pair $p,q$ satisfying
\begin{equation} \label{TH1pq}
  1\le 3/p+2/q \le 2, \quad 2<q \le \infty, \quad (p,q) \not = (3/2, \infty),
\end{equation}
there exists a constant $\eps>0$ depending only on $p,q,\gamma$
and $\norm{f}_{M_{2,\gamma}}$ such that, if the pair $u,p$ is a
suitable weak solution of the Navier-Stokes equations
\eqref{intro:nse} vanishing on a flat boundary $\Gamma$ according to
Definition \ref{defsws:100}, and for some point
$z=(x,t)\in\Gamma\times (0,T)$, $u$ is locally in $L^{p,q}$ near
$z$ and
\begin{equation}
\limsup_{r\raro 0_+}
r^{-( \frac 3p + \frac 2q -1)}
\norm{ \norm{u(y,s)}_{L^p(B^+_{x,r})} }_{L^q(t-r^2,t)}
\leq \eps, 
\label{nse:abs30}
\end{equation}
then $z$ is a regular point.
\end{thm}

{\bf Comments for Theorem \ref{TH1}}:

\begin{enumerate}
\item The same statement for an interior point $z$
remains true, see Appendix.

\item The quantities in \eqref{nse:abs30} are invariant under the
scaling $u(x,t) \to s u( sx,s^2t)$. Scaling invariant quantities have
been important in the study of \eqref{intro:nse}, see e.g. \cite{CKN}.

\item The exponents $(p,q)$ in Theorem \ref{TH1} correspond to Region
II in Figure 1, which is a solid parallelogram excluding its top
borderline and the corner point $(2/3,0)$.  By H\"older's inequality, it
suffices to prove the cases $\frac 3p + \frac 2q = 2$, $2<q<\infty$,
the right borderline of Region II.
Our method fails for the end points $(p,q)=(3/2,\infty)$ and $(3,2)$ for the
lack of $L^{3/2,1}$ and $L^{1,2}$ estimates for the Stokes system.

\item The usual \ps conditions correspond to Region I and imply
\eqref{nse:abs30} pointwise. Thus, also by H\"older's inequality,
regularity under the \ps conditions is a corollary of Theorem \ref{TH1},
except in the endpoint cases $(p,q)=(3,\infty)$ or $(p,q)=(\infty,2)$ 
(see Corollary \ref{swecor} for the details).
Regularity up to the boundary under the \ps conditions is proved in
\cite{k_bnse,VS02} but the proof of Theorem \ref{TH1} seems easier.

\item One key feature of Theorem \ref{TH1} is that {\bf condition
\eqref{nse:abs30} does not involve any scaled norm of the pressure
$p$}.  A previous such result is by Tian and Xin \cite{TX} for the
special case of \eqref{nse:abs30} with $(p,q)=(3,3)$.  Another such
result is by Seregin and \v Sver\'ak \cite{SS} for
$(p,q)=(2,\infty)$. Both results are for interior points,
and are included in Region II.

\item Another regularity criterion in \cite{TX} is the uniform
boundedness \\
$
\sup_{r < R_0} \bke{r^{-1/2} \norm{u}_{L^{2,\infty}(Q_{z,r})} }\le M
$
for some $M< \infty$, and the condition \eqref{nse:abs30} with
$(p,q)=(2,2)$ and a small constant $\e$ depending on $M$. Although
$(p,q)=(2,2)$ lies outside of Region II, using $\displaystyle
\norm{u}_{L^{2,4}} \le \norm{u}_{L^{2,\infty}}^{1/2}
\norm{u}_{L^{2,2}} ^{1/2}$,
one obtains \eqref{nse:abs30} with $(p,q)=(2,4)$, which
falls in Region II. Thus this result is also implied by Theorem
\ref{TH1}.

\item
Eq.~\eqref{nse:abs30} is a uniform estimate for $r$ sufficiently
small.  There are conditions which only require one $r$. For
example, there is an $\e>0$ such that the condition
\[
  r^{-2} \int_{Q_{z,r}} \bke{|u|^3 + |p|^{3/2}} dx\,dt \le \e \quad
  \text{ for some } r>0
\]
implies regularity at $z$. This is essentially \cite[Proposition
1]{CKN} and is stated as above in \cite{NRS96,L98}. Also see
\cite{S02} and our Lemma \ref{mod:lemma} when $z$ is on boundary.

\end{enumerate}


{
\setlength{\unitlength}{2mm}
\noindent
\begin{picture}(30,28)
\put (1,0){\vector(1,0){26}}
\put (1,0){\vector(0,1){26}}
\put (1,24){\line(1,0){24}}
\put (25,0){\line(0,1){24}}
\put (9,0){\line(-2,3){8}}
\put (17,0){\line(-2,3){8}}
\put (1,12){\line(1,0){8}}
\put(28,0){\makebox(0,0)[c]{\scriptsize $\frac 1p$}}
\put(0,26){\makebox(0,0)[c]{\scriptsize $\frac 1q$}}
\put(1,-1){\makebox(0,0)[c]{\scriptsize $0$}}
\put(9,-1){\makebox(0,0)[c]{\scriptsize $1/3$}}
\put(11.5,12){\makebox(0,0)[c]{\scriptsize $(\frac 13, \frac 12)$}}
\put(17,0){\makebox(0,0)[c]{\scriptsize o}}
\put(9,12){\makebox(0,0)[c]{\scriptsize o}}
\put(9,8){\makebox(0,0)[c]{\scriptsize $\bullet$}}
\put(15,8){\vector(-1,0){5}}
\put(17.5,8){\makebox(0,0)[c]{\scriptsize $(\frac 13, \frac 13)$}}
\put(13,-1){\makebox(0,0)[c]{\scriptsize $1/2$}}
\put(13,0){\makebox(0,0)[c]{\scriptsize $\bullet$}}
\put(17,-1){\makebox(0,0)[c]{\scriptsize $2/3$}}
\put(25,-1){\makebox(0,0)[c]{\scriptsize $1$}}
\put(0,0){\makebox(0,0)[c]{\scriptsize $0$}}
\put(0,12){\makebox(0,0)[c]{\scriptsize $\frac 12$}}
\put(0,24){\makebox(0,0)[c]{\scriptsize $1$}}
\put(3,4){\makebox(0,0)[c]{\scriptsize I}}
\put(10,4){\makebox(0,0)[c]{\scriptsize II}}
\put(13,-3){\makebox(0,0)[c]{\scriptsize Figure 1: Regularity Criteria}}
\end{picture}
}
\quad
{
\setlength{\unitlength}{2mm}
\noindent
\begin{picture}(30,28)
\put (1,0){\vector(1,0){26}}
\put (1,0){\vector(0,1){26}}
\put (1,24){\line(1,0){24}}
\put (25,0){\line(0,1){24}}
\put (9,0){\line(-2,3){8}}
\put (13,0){\line(-2,3){8}}
\put (1,12){\line(1,0){4}}
\put (1,12){\line(1,-1){6}}
\put (9,0){\line(-1,3){2}}
\put (5,12){\line(-1,3){4}}
\put(13,6){\vector(-1,0){5}}
\put(15.5,6){\makebox(0,0)[c]{\scriptsize $(\frac 14, \frac 14)$}}
\put(7,6){\makebox(0,0)[c]{\scriptsize o}}
\put(28,0){\makebox(0,0)[c]{\scriptsize $\frac 1p$}}
\put(0,26){\makebox(0,0)[c]{\scriptsize $\frac 1q$}}
\put(1,-1){\makebox(0,0)[c]{\scriptsize $0$}}
\put(9,-1){\makebox(0,0)[c]{\scriptsize $1/3$}}
\put(13,-1){\makebox(0,0)[c]{\scriptsize $1/2$}}
\put(25,-1){\makebox(0,0)[c]{\scriptsize $1$}}
\put(0,0){\makebox(0,0)[c]{\scriptsize $0$}}
\put(0,12){\makebox(0,0)[c]{\scriptsize $\frac 12$}}
\put(0,24){\makebox(0,0)[c]{\scriptsize $1$}}
\put(5,12){\makebox(0,0)[c]{\scriptsize $\bullet$}}
\put(7.5,12){\makebox(0,0)[c]{\scriptsize $(\frac 16, \frac 12)$}}
\put(13,0){\makebox(0,0)[c]{\scriptsize $\bullet$}}
\put(1,24){\makebox(0,0)[c]{\scriptsize o}}
\put(3,4){\makebox(0,0)[c]{\scriptsize I}}
\put(2.5,14){\makebox(0,0)[c]{\scriptsize III}}
\put(15,14){\makebox(0,0)[c]{\scriptsize IV}}
\put(7,4){\makebox(0,0)[c]{\scriptsize V}}
\put(9,4){\makebox(0,0)[c]{\scriptsize VI}}
\put(13,-3){\makebox(0,0)[c]{\scriptsize Figure 2: Partial Regularity}}
\end{picture}
}

\vspace{10mm}


The main tools of our analysis are a standard ``blow up'' method and
the decomposition of the pressure as introduced in \cite{S02}, which
enable us to prove a decay property of the scaled Lebesgue norms of
velocity and pressure in both the interior and boundary cases (see
Lemma \ref{lemma:seregin} and Appendix).  Combining this with 
the local estimate
of the Stokes system for the pressure, we can estimate the pressure
for the Navier-Stokes equations near the boundary (see Lemma
\ref{preest:100}).

As mentioned earlier, the best available estimate for the singular set
is that $\calP^1(\calS)=0$ (in \cite{CL00} the estimate of the Hausdorff 
measure of the singular set for suitable weak solutions
was improved by a logarithmic factor for the interior case). 
In the following theorem we improve the
estimate using Theorem \ref{TH1}, assuming some conditions of the \ps
type.

\begin{thm}[Partial regularity] \label{TH2}
Suppose $f\in M_{2,\gamma}(Q)$ for some $\gamma>0$, a parabolic
Morrey space.  Suppose $(u,p)$ is a weak solution of the
Navier-Stokes equations \eqref{intro:nse} according to Definition 2.2 and
assume that
\begin{equation*}
u\in L^m((0,T);L^l(\om)),\quad
(l,m) \in \text{Region V},
\end{equation*}
where Region V is the triangular region in Figure 2 satisfying $
\frac{3}{l}+\frac{2}{m}>1$, $\frac{1}{l}+\frac{1}{m} < \frac 12$, and
$\frac{3}{l}+\frac{1}{m}<1$.  Let $S$ denote the singular set of $u$
up to the flat boundary where $u$ vanishes, and
\begin{equation*}
d(l,m) = \left \{
\begin{aligned}
&3-m+ \frac {2m}l \quad & \text{if }\ l > m ,\\
&2 -m + \frac {3m} l
& \text{if }\  l \leq m.
\end{aligned} 
\right .
\end{equation*}
Then the $d(l,m)$ dimensional parabolic Hausdorff measure of $S$ is zero,
$\calP ^{d(l,m)}(S)=0$.
\end{thm}

In Theorem \ref{TH2} we only require a weak solution. As for the \ps class,
these solutions are automatically suitable weak solutions, see section 4.
Note that weak solutions are known to lie in Region IV of Figure 2,
including the solid line from $(\frac 12,0)$ to $(\frac 16,\frac 12)$.
Assuming $u\in L^{p,q}$, one can use Theorem \ref{TH1} to estimate the
dimension of the singular set for all $(p,q)$ in Region
II with $q < \infty$, 
but only those in Region V give us dimensions less than 1.

The plan of this paper is as follows:
In Section 2 we introduce the notion of suitable weak solutions near the
boundary, which is a slightly modified version of that used in \cite{S02}.
We also show the decay property of the velocity field and pressure 
(see Lemma \ref{mod:lemma} and Lemma \ref{lemma:seregin}).  
In Section 3 we present the proof of the main Theorem \ref{TH1}.
In Section 4, as an application, we investigate the size of the possible 
singular set under our additional integrability assumption on $u$
(see Theorem \ref{TH2}).
In the Appendix we present a brief sketch of the proof that 
the regularity criteria \eqref{nse:abs30} is valid in the interior.


\section{Preliminaries}
In this section we introduce notation, define suitable weak
solutions, and give some lemmas on the decay properties of the
velocity and pressure.

We start with notation.  Denote by $\om$ an open domain in
$\R^3$ and by $\pr\om$ its boundary. $\Gamma$ indicates an open
subset of $\pr\om$ which lies on a plane.  In this article,
for simplicity, we assume $\Gamma$ lies on the plane $\{x_3=0\}$.

For $1\leq q\leq\infty$, $W^{k,q}(\om)$ denotes the usual Sobolev
space, i.e. $W^{k,q}(\om)= \{ u\in L^q(\om):D^{\ap}u\in L^q(\om),
0\leq|\ap|\leq k\}$.  As usual, $W^{k,q}_0(\om)$ is the completion of
$\calC^{\infty}_0(\om)$ in the $W^{k,q}(\om)$ norm.  We also denote by
$W^{-k,q'}(\om)$ the dual space of $W^{k,q}_0(\om)$ where $q$ and $q'$
are H\"older conjugates.

For a domain $Q\subset\R^3\times I$, we denote by
$\calC^{\ap,\frac{\ap}{2}}_{x,t}(Q)$ the Banach space of functions
that are H\"older continuous with exponent $\ap\in (0,1)$, with
respect to the parabolic metric $d(z,z')=|x-x'|+|t-t'|^{\frac{1}{2}}$
where $z=(x,t)$ and $z'=(x',t')$.

We denote by $\aint_E f$ the average of $f$ on $E$; i.e.,
$\aint_E f=\int_E f/|E|$.

We denote by $M_{2,\gamma}$ a parabolic version of Morrey's spaces
(see e.g. \cite[page 3]{S02}). For $\omega\subset\R^3\times\R$ 
and a positive number $\gamma \in (0,2]$, we define the space
\begin{eqnarray*}
M_{2,\gamma}(\omega;\R^3):=\{f\in L_{2,\rm{loc}}(\omega;\R^3):
m_{\gamma}(f;\omega)<\infty\},
\end{eqnarray*}
where 
\be
m_{\gamma}(f;\omega):=\sup\bket{\frac{1}{r^{\gamma-2}}
\bke{\aint_{Q(z,r)\cap\omega}\abs{f}^2dz'}^{\frac{1}{2}}:z\in\bar{\omega}, r>0}.
\label{morrey:100}
\ee

Parabolic Hausdorff measures are defined in \cite{CKN} using parabolic
cylinders instead of usual balls. For any $X \subset \R^3 \times \R$
and $k \ge 0$ one defines
\begin{equation*}
\calP^k(X) = \inf_{\delta \to 0^+} \calP^k_\delta(X),
\quad
\calP^k_\delta(X) = \inf \left \{\sum_{i = 1}^\infty 
r_i^k: X \subset \bigcup _i Q_{z_i,r_i}, \ r_i < \delta \right \}.  
\end{equation*}

Finally, by $N=N(\alpha,\beta,\ldots)$ we denote a constant depending
on the prescribed quantities $\alpha,\beta,\ldots$, which may change
from line to line.

\medskip

Next, we define several scaling-invariant functionals.  Let
$z=(x,t)\in\Gamma\times I$.  As in \cite{CKN, L98, LS99, S02}, let
\begin{equation*}
A(r):=\sup_{t-r^2\leq s<t}
\frac{1}{r}\int_{B^+_{x,r}}\abs{u(y,s)}^2dy,
\end{equation*}
\begin{equation*}
C(r):=\frac{1}{r^2}\int_{Q^+_{z,r}}\abs{u(y,s)}^3\ind y\ind s,
\end{equation*}
\begin{equation*}
E(r):=\frac{1}{r}\int_{Q^+_{z,r}}\abs{\nb u(y,s)}^2\ind y\ind s.
\end{equation*}
Let $\kappa, \kappa^*$ and $\lambda$ be numbers satisfying
\be
\frac{3}{\kappa}+\frac{2}{\lambda}=4,\qquad 
\frac{1}{\kappa^*}=\frac{1}{\kappa}-\frac{1}{3},\qquad
1<\lambda<2.
\label{kalam:10}
\ee
We also introduce new functionals, 
which are useful for us:
\be
\tilde{D}(r):= \frac{1}{r}
\bke{\int_{t-r^2}^{t}
\bke{\int_{B^+_{x,r}}\abs{p(y,s)-(p)_{B^+_{x,r}}(s)}^{\kappa^*}
\ind y}^{\frac{\lambda}{\kappa^*}}
\ind s}^{\frac{1}{\lambda}},
\label{tild:p10}
\ee
where $(p)_{B^+_{x,r}}(s)=\aint_{B^+_{x,r}}p(y,s)\ind y$, 
\be
\tilde{D}_1(r):= \frac{1}{r}
\bke{\int_{t-r^2}^{t}
\bke{\int_{B^+_{x,r}}\abs{\nb p(y,s)}^{\kappa}\ind y}^{\frac{\lambda}{\kappa}}
\ind s}^{\frac{1}{\lambda}},
\label{tild:p20}
\ee
and finally,
\begin{equation*}
G(r):= \frac{1}{r}\bke{\int_{t-r^2}^t
\bke{\int_{B^+_{x,r}}\abs{u(y,s)}^p\ind y}^{\frac{q}{p}}\ind s}^{\frac{1}{q}},
\end{equation*}
where $p$ and $q$ are the H\"older conjugate exponents of $\kappa^*$
and $\lambda$ in \eqref{kalam:10}, i.e.  \be
\frac{1}{p}+\frac{1}{\kappa^*}=1, \qquad
\frac{1}{q}+\frac{1}{\lambda}=1.
\label{kalampq:100}
\ee
It is straightforward from \eqref{kalam:10} that $p$ and $q$ satisfy 
\be
\frac{3}{p}+\frac{2}{q}=2,\qquad 2<q<\infty.
\label{reg:20}
\ee

\begin{rem}
In \cite{S02} the following functionals are used, instead of 
$\tilde{D}(r), \tilde{D}_1(r)$,
\[
D(r):=\frac{1}{r^2}\int_{Q^+_{z,r}}\abs{p-(p)_{B^+_{x,r}}}^{\frac{3}{2}}
\ind z,
\quad
D_1(r):= \frac{1}{r^{\frac{3}{2}}}\int_{t-r^2}^{t}
\bke{\int_{B^+_{x,r}}\abs{\nb p}^{\frac{9}{8}}\ind y}^{\frac{4}{3}}\ind s.
\]
We note that $D_1(r)$ is a special case of $\tilde{D}_1$ with
$\kappa=\frac{9}{8}, \lambda=\frac{3}{2}$.  \qed
\end{rem}

Next we define suitable weak solutions for the Navier-Stokes equations.
\begin{defi}\label{defsws:100}
Let $Q=\om\times I$ where $\om\subset\R^3$ and $I=[0,T)$ and $\Gamma$
be an open subset of the set $\pr\om$.  Suppose that $f$ belongs to
the parabolic Morrey space $M_{2,\gamma}(Q)$ for some $\gamma \in (0,2]$.  A
pair of ($u$, $p$) is a {\bf suitable weak solution} to the
Navier-Stokes equation (\ref{intro:nse}) in $Q$ near the boundary
$\Gamma$ and vanishing on $\Gamma$ if the following conditions are
satisfied.
\begin{itemize}
\item[(a)]
The functions $u:Q \to \R^3$ and $p:Q \to \R$ satisfy
\bes
u\in L^{\infty}(I;L^2(\om))\cap L^2(I;W^{1,2}(\om)),
\ees
\be
p\in L^{\lambda}(I;L^{\kappa^*}(\om)),
\label{defsws:200}
\ee
\be
\nb^2 u\in L^{\lambda}(I;L^{\kappa}(\om)),\quad
\nb p\in L^{\lambda}(I;L^{\kappa}(\om)),
\label{defsws:300}
\ee
where $\kappa, \kappa^*$ and $\lambda$ are fixed numbers 
satisfying \eqref{kalam:10}. 

\item[(b)]
$u$ and $p$ solve the Navier-Stokes equations (\ref{intro:nse}) in $Q$ 
in the sense of distributions and $u$ satisfies the boundary condition
$u=0$ on $\Gamma\times I$.

\item[(c)]
$u$ and $p$ satisfy the local energy inequality

\parbox{11cm}{
\[
\int_{\om}\abs{u(x,t)}^2\phi(x,t)\ind x
+2\int_{Q}\abs{\nb u(x,t')}^2\phi(x,t')\ind x\ind t'
\]
\[
\leq\int_{Q}\bke{\abs{u}^2(\pr_t\phi+\Delta\phi)+(\abs{u}^2+2p)u\cdot\nb\phi
+2f\cdot u\phi}\ind x\ind t' 
\]} \hfill
\parbox{1cm}{\begin{equation}\label{lei}\end{equation}} 
for almost all $t\in (0,T)$ and for all nonnegative functions 
$\phi\in\calC^{\infty}_0(\R^3 \times \R)$, 
vanishing in a neighbourhood of the set
$\bke{\om\times \{ t=0 \}}\cup \bke{(\pr\om\setminus\Gamma)\times (0,T)}$.
\end{itemize}
\end{defi}

Let us make several comments on the above definition.

\begin{rem}
Sohr and Von Wahl \cite{SW86} showed that, under reasonable
assumptions on $f$ and $u_0$, the pressure $p$ of a weak solution
belongs to $L^{\frac{5}{3}}(\om \times I)$, which corresponds to
$\kappa^*=\lambda=\frac 53$ in \eqref{defsws:200}.  Here
$\om\subset\R^3$ can be either a bounded domain, an exterior domain,
or a half-space.  Giga and Sohr \cite{GS} later proved that $u_t,
\nb^2 u, \nb p\in L^{\kappa, \lambda}(Q)$ and $p\in L^{\kappa^*,
\lambda}(Q)$ where $\kappa, \kappa^*$ and $\lambda$ are any numbers
satisfying \eqref{kalam:10}.  Therefore, it seems reasonable to make
assumptions \eqref{defsws:200} and \eqref{defsws:300} for suitable
weak solutions.
\end{rem}

\begin{rem}
The main difference between suitable weak solutions and the original
Leray-Hopf weak solutions is the additional condition of the local
energy inequality \eqref{lei}.  The existence of suitable weak
solutions is proved in \cite{VS77,CKN}.  Slightly modified
definitions are used in \cite{L98,LS99,S02}.  As indicated in
\cite[Remarks 4, page 823]{CKN}, it is an open question if all weak
solutions are suitable weak solutions.
\end{rem}

Next we show the local regularity criterion near the boundary, which is 
analogous to Proposition 2.6 in \cite{S02}. Although 
our proof is based on a standard ``blow up''
method similar to that of \cite{S02}, 
we present its details since different functionals are used
for the pressure, and therefore modifications are needed.

\begin{lemm}\label{mod:lemma}
There exists $\eps > 0$ depending only on 
$\lambda$, $\gamma$ and
$\norm{f}_{M_{2,\gamma}}$, such that if $u$ is a 
suitable weak solution of the Navier-Stokes equations 
satisfying Definition \ref{defsws:100}, $z=(x,t)\in\Gamma\times I$, and
\[
  \liminf_{r\raro 0_+}\bke{C^{\frac{1}{3}}(r)+\tilde{D}(r)}<\eps,
\] 
then $z$ is a regular point.
\end{lemm}
Before we prove Lemma \ref{mod:lemma}, 
we first prove the following lemma, which gives a decay property of
$u$ and $p$ in some Lebesgue spaces.

\begin{lemm}\label{lemma:seregin}
Let $0<\theta<1/2$ and $0 < \beta < \gamma \le 2$.
There exist $\eps_1, r_1>0$ depending on $\lambda,\theta, \gamma$ and $ \beta$
such that if $u$ is a suitable weak solution of the Navier-Stokes equations 
satisfying Definition \ref{defsws:100}, $z=(x,t)\in\Gamma\times I$,
and $C^{\frac{1}{3}}(r)+\tilde{D}(r)+ m_{\gamma}(f) r^{\beta+1}
<\eps_1$ for some $r \in (0, r_1)$,
then
\[
\bke{C^{\frac{1}{3}}(\theta r)+\tilde{D}(\theta r)}
<N\theta^{1+\ap}\bke{C^{\frac{1}{3}}(r)+\tilde{D}(r) 
+  m_{\gamma}(f) r^{\beta+1}},
\]
where $0<\ap<1$ and $N>0$ are absolute constants.
\end{lemm}
\begin{pf}
 For simplicity 
we assume $f=0$. The general case follows similarly.
For convenience, we denote $\varphi(r):=C^{\frac{1}{3}}(r)+\tilde{D}(r)$.
Suppose the statement is not true. 
Then for any $\alpha \in (0,1)$ and $N > 0$,
there exist $z_n=(x_n,t_n)$, $r_n\searrow 0$, and $\eps_n\searrow 0$ 
such that 
\be
\varphi(r_n)=\eps_n,\quad \varphi(\theta r_n)
>N\theta^{1+\ap}\varphi(r_n)=N\theta^{1+\ap}\eps_n.
\label{13feb100}
\ee
Let $w=(y,s)$ where $y=r^{-1}_n(x-x_n), s=r^{-2}_n (t-t_n)$ 
and we define $v_n$ and $q_n$ as follows:
\[
v_n(w)=\eps^{-1}_n r_n u(z),\quad 
q_n(w)=\eps^{-1}_n r^2_n \bke{p(z)-(p)_{B^+_{r_n}}(t)}.
\]
For convenience we also define $C(v_n, \theta)$, 
$\tilde{D}(q_n, \theta)$, and  $\tilde{D}_1(q_n, \theta)$ by
\[
C(v_n, \theta):=\frac{1}{\theta^2}\int_{Q^+_{\theta}}\abs{v_n}^3 dw,\quad
\tilde{D}(q_n, \theta):= \frac{1}{\theta}
\bke{\int_{-\theta^2}^{0}\bke{\int_{B^+_{\theta}}
\abs{q_n-(q_n)_{B^+_{\theta}}}^{\kappa^*}dy}^{\frac{\lambda}{\kappa^*}}
ds}^{\frac{1}{\lambda}},
\]
\[
\tilde{D}_1(q_n, \theta):=\frac{1}{\theta}
\bke{\int_{-\theta^2}^{0}\bke{\int_{B^+_{\theta}}
\abs{\nb q_n}^{\kappa}dy}^{\frac{\lambda}{\kappa}}
ds}^{\frac{1}{\lambda}},
\]
where $\kappa^*, \kappa$ and $\lambda$ are numbers in \eqref{kalam:10}.
By the change of variables, we have
\be
\frac{1}{\eps_n}\varphi(\theta r_n)
=C^{\frac{1}{3}}(v_n, \theta)+\tilde{D}(q_n, \theta).
\label{13feb200}
\ee
For convenience, we denote 
$\psi_n(\theta):=C^{\frac{1}{3}}(v_n, \theta)+\tilde{D}(q_n, \theta)$.
Due to \eqref{13feb100} and \eqref{13feb200},  we get
\be
\psi_n(1)=\norm{v_n}_{L^3(Q^+_1)}+\norm{q_n}_{L^{\kappa^*,\lambda}(Q^+_1)}=1,
\label{weakw:100}
\ee
\be
\psi_n(\theta)=C^{\frac{1}{3}}(v_n, \theta)+\tilde{D}(q_n, \theta)
\geq N\theta^{1+\ap}.
\label{weakq:100}
\ee
On the other hand, $v_n, q_n$ solve the following system in a weak sense
\[
\pr_s v_n-\Delta v_n+\eps_n (v_n \cdot \nb)v_n+\nb q_n=0,\quad
\rm{div}\,\,v_n=0\qquad \mbox{ in }\,\,\, Q^+_1
\]
with
\[
v_n=0 \quad\mbox{ on }\,\,\,(B_1\cap\{x_3=0\})\times (-1,0).
\]
Because of \eqref{weakw:100}, we have following weak convergence (possibly 
subsequences of $v_n$ and $q_n$ should be taken, however
we use the same symbol for simplicity)  
\[
v_n\raroup v \mbox{ in }L^3(Q^+_1),\quad 
q_n\raroup q \mbox{ in }L^{\kappa^*,\lambda}(Q^+_1),\quad
(q)_{B^+_1}(s)=0.
\]
In addition, one can easily see that 
$\pr_s v_n$ is uniformly bounded
in $L^{\lambda}\bke{(-1,0);(W^{2,2}(B^+_1))'}$, 
and, therefore, we also have
\begin{equation}
\label{eq:weak}
\pr_sv_n\raroup \pr_sv \,\,\mbox{ in }\,\,
L^{\lambda}\bke{(-1,0);\bkp{W^{2,2}(B^+_1)}'}. 
\end{equation}
Next we show that $\nb v_n$ is uniformly bounded in 
$L^{2}(Q^+_{\frac{3}{4}})$.
Let 
$\phi$ be a standard cut off function satisfying $\phi$ is smooth,
\[
\phi=1\,\,\,\mbox{ on }Q_{\frac{3}{4}},
\quad\phi=0 \,\,\,\mbox{ on }
\bke{\R^3\times(-\infty,0)}\setminus Q_1,\quad 0\leq\phi\leq 1.
\]
From the local energy inequality, for every $\tau\in (-1,0)$,
we obtain
\[
\int_{B^+_1}\abs{v_n(\cdot,\tau)}^2\phi^2(x,\tau)\ind y
+\int_{-1}^{\tau}
\int_{B^+_1}\abs{\nb v_n}^2\phi^2\ind y\ind s
\]
\[
\leq N\Big(\int_{-1}^{\tau}\int_{B^+_1}\abs{v_n}^2
\bke{\abs{\pr_s\phi}+\abs{\Delta\phi}+\abs{\nb\phi}}\ind y\ind s+
\]
\[
+\eps_n\int_{-1}^{\tau}\int_{B^+_1}\abs{v_n}^3
\abs{\nb\phi}\ind y\ind s
+\int_{-1}^{\tau}\int_{B^+_1}\abs{q_nv_n\cdot\nb\phi\phi}\ind y\ind s\Big).
\]
Consider the last term in the above inequality. Using the H\"older 
inequality, we have
\[
\int_{-1}^{\tau}\int_{B^+_1}\abs{q_nv_n\cdot\nb\phi\phi}\leq 
\bke{\int_{-1}^{\tau}\big(\int_{B^+_1}\abs{q_n\nb\phi}^{\kappa^*}\big)
^{\frac{\lambda}{\kappa^*}}}^{\frac{1}{\lambda}}
\bke{\int_{-1}^{\tau}\big(\int_{B^+_1}\abs{v_n\phi}^p
\big)^{\frac{q}{p}}}^{\frac{1}{q}}.
\]
where $\kappa^*, \lambda, p$, and $q$ are numbers in 
\eqref{kalam:10}, \eqref{kalampq:100}, and \eqref{reg:20}. 
We recall that $p$ and $q$ are in the ranges $3/2<p<3$ and $2<q<\infty$.
In case $q\leq 3$, since $p,q\leq 3$, we have
\[
\norm{v_n\phi}_{L^{p,q}\bkp{(-1,\tau);B^+_1}}
\leq N\norm{v_n\phi}_{L^3(Q^+_1)}.
\]
Therefore, in this case, $\nb v_n$ is uniformly bounded in 
$L^2(Q^+_{3/4})$
because of \eqref{weakw:100}.
It remains to consider the case $3<q<\infty$ (equivalently $3/2<p<9/4$).
Suppose $2<p<9/4$, which is equivalent to $4<q<\infty$. 
In this case, by interpolation, one can see
the following estimate:
\[
\norm{v_n\phi}_{L^{p,q}(Q^+_1)}\leq N\sup_{-1<s<\tau}
\norm{v_n\phi(\cdot,s)}^{\frac{2\ap}{p}}_{L^{2}(B^+_1)}
\norm{v_n\phi}^{\frac{3(1-\ap)}{p}}_{L^3(Q^+_1)},
\]
where $\ap=3-p$. In the above inequality, we used that $(1-\ap)q/p<1$.
Since $\tau$ is arbitrary between $-1$ and $0$, 
we obtain
\[
\sup_{-1<\tau<0}\norm{v_n\phi(\cdot,\tau)}_{L^{2}(B^+_1)}
+\norm{\nb v_n\phi}_{L^{2}(Q^+_1)}\leq N\Big(\norm{v_n}_{L^{2}(Q^+_1)}+
\]
\[
+\eps_n\norm{v_n}_{L^{3}(Q^+_1)} 
+\sup_{-1<s<0}\norm{v_n\phi(\cdot,s)}^{\frac{2(3-p)}{p}}_{L^{2}(B^+_1)}
\norm{v_n\phi}^{\frac{3(p-2)}{p}}_{L^3(Q^+_1)}
\norm{q_n}_{L^{\kappa^*,\lambda}(Q^+_1)}\Big).
\]
Therefore, using Young's inequality, we obtain
\[
\sup_{-1<\tau<0}\norm{v_n\phi(\cdot,\tau)}_{L^{2}(B^+_1)}
+\norm{\nb v_n\phi}_{L^{2}(Q^+_1)}\leq N\Big(\norm{v_n}_{L^{2}(Q^+_1)}+
\]
\[
+\eps_n\norm{v_n}_{L^{3}(Q^+_1)}+\norm{v_n\phi}_{L^3(Q^+_1)}
\norm{q_n}^{\frac{p}{3(p-2)}}_{L^{\kappa^*,\lambda}(Q^+_1)}
\Big).
\]
Therefore, we have 
\[
\norm{\nb v_n}_{L^{2}(Q^+_{\frac{3}{4}})}\leq 
N\Big(\norm{v_n}_{L^{2}(Q^+_1)}+\eps_n\norm{v_n}_{L^{3}(Q^+_1)} +
\norm{v_n\phi}_{L^3(Q^+_1)}
\norm{q_n}^{\frac{p}{3(p-2)}}_{L^{\kappa^*,\lambda}(Q^+_1)}\big).
\]
Therefore, $\nb v_n$ is also uniformly bounded in 
$L^{2}(Q^+_{\frac{3}{4}})$ for the case $2<p<9/4$.
For the case $3/2<p\leq 2$ (equivalently $3<q\leq 4$) we have
\[
\norm{v_n\phi}_{L^{p,q}(Q^+_1)}\leq N\sup_{-1<s<\tau}
\norm{v_n\phi(\cdot,s)}_{L^{2}(B^+_1)}.
\]
By proceeding as for the previous case, we can 
obtain the uniform bound of $\nb v_n$ in 
$L^{2}(Q^+_{\frac{3}{4}})$.
So together with~(\ref{eq:weak}), we get 
\[
\nb v_n\,\,\raroup\,\, \nb v\,\,\,\mbox{ in }L^2(Q^+_{\frac{3}{4}}),\qquad
v_n\,\,\raro\,\, v \,\,\,\mbox{ in }L^3(Q^+_{\frac{3}{4}}).
\]
Moreover, $v$ and $q$ solve the following linear Stokes system
\[
\pr_s v-\Delta v+\nb q=0,\quad
{\rm{div}}\, v=0\qquad \mbox{ in }\,\,\, Q^+_1
\]
with
\[
v=0 \quad\mbox{ on }\,\,\,(B_1\cap\{x_3=0\})\times (-1,0).
\]
Next we show that 
\be
\pr_s v_n, \nb^2 v_n, \nb q_n\raroup \pr_s v, \nb^2 v, \nb q
\quad \mbox{ in }L^{\kappa,\lambda}(Q^+_{\frac{5}{8}})\mbox{, respectively}.
\label{wq:100}
\ee
Indeed, after direct calculations, we obtain
\be
\norm{(v_n\cdot\nb)v_n}_{L^{\kappa,\lambda}(Q^+_{\frac{3}{4}})}
\leq N\norm{\nb v_n}^{\frac{2}{\lambda}}_{L^2(Q^+_{\frac{3}{4}})}
\norm{v_n}^{\frac{3-2\kappa}{\kappa}}_{L^{2,\infty}(Q^+_{\frac{3}{4}})}.
\label{preest:890}
\ee
Due to the local boundary estimate for the Stokes system 
(see \cite[Proposition 1]{S00}), 
we have the following estimate for $v_n$ and $q_n$:
\[
\norm{\pr_s v_n}_{L^{\kappa,\lambda}(Q^+_{\frac{5}{8}})}
+\norm{\nb^2 v_n}_{L^{\kappa,\lambda}(Q^+_{\frac{5}{8}})}
+\norm{\nb q_n}_{L^{\kappa,\lambda}(Q^+_{\frac{5}{8}})}
\]
\[
\leq N\bke{ \norm{v_n}_{L^{\kappa,\lambda}(Q^+_{\frac{3}{4}})}
+\norm{\nb v_n}_{L^{\kappa,\lambda}(Q^+_{\frac{3}{4}})}
+\norm{q_n}_{L^{\kappa,\lambda}(Q^+_{\frac{3}{4}})}+
\eps_n\norm{(v_n\cdot\nb)v_n}_{L^{\kappa,\lambda}(Q^+_{\frac{3}{4}})}}.
\]
Therefore, we obtain
\be
\norm{\pr_s v_n}_{L^{\kappa,\lambda}(Q^+_{\frac{5}{8}})}
+\norm{\nb^2 v_n}_{L^{\kappa,\lambda}(Q^+_{\frac{5}{8}})}
+\norm{\nb q_n}_{L^{\kappa,\lambda}(Q^+_{\frac{5}{8}})}
\leq N(1+\eps_n),
\label{preest:900}
\ee
where we used \eqref{preest:890}.
The assertion \eqref{wq:100} is established.

According to H\"older estimate of the 
Stokes system near boundary (see \cite[Lemma 1]{S00}),
$v$ is H\"older continuous in $Q^+_{\frac{1}{2}}$ with the exponent $\ap$
with $0<\ap<2(1-1/\lambda)$. Here we fix $\ap=1-1/\lambda$, 
denoted by $\ap_0$ from now on.
Then, by H\"older continuity of $v$ and strong convergence
of the $L^3$-norm of $v_n$, we obtain
\be
C^{\frac{1}{3}}(v,\theta)\leq N_1\theta^{1+\ap_0},\quad 
C(v_n,\theta)\,\,\raro\,\, C(v,\theta).
\label{hs:100}
\ee
Let $\tilde{B}^+$ be a domain with smooth boundary such that 
$B^+_{11/16}\subset\tilde{B}^+\subset B^+_{3/4}$, and 
$\tilde{Q}^+:=\tilde{B}^+\times (-(3/4)^2,0)$.
Now we consider the following initial and boundary problem:
\[
\pr_s \hat{v}_n-\Delta\hat{v}_n +\nb\hat{q}_n
=-\eps_n (v_n \cdot \nb)v_n,\quad
\rm{div}\,\,\hat{v}_n=0 \quad \mbox{ in }\tilde{Q}^+,
\]
\[
(\hat{q}_n)_{\tilde{B}^+}(s)=0, \quad s\in 
\bke{-\bkp{{\frac{3}{4}}}^2,0},
\]
\[
\hat{v}_n=0 \,\,\,\mbox{ on }\pr\tilde{B}^+
\times \bkt{-({\frac{3}{4}})^2,0},\quad 
\hat{v}_n=0 \,\,\mbox{ on }\,\,\tilde{B}^+\times
\bket{s=-({\frac{3}{4}})^2}
\]
Using the global estimate of the Stokes system (see \cite[Theorem 3.1]{GS}), 
we get the following estimate 
\[
\norm{\pr_s\hat{v}_n}_{L^{\kappa,\lambda}(\tilde{Q}^+)}
+\norm{\hat{v}_n}_{L^{\kappa}\bkp{(-({\frac{3}{4}})^2,0);
W^{2,\lambda}_0(\tilde{B}^+)}}
+\norm{\hat{q}_n}_{L^{\kappa}\bkp{(-({\frac{3}{4}})^2,0);
W^{1,\lambda}(\tilde{B}^+)}}
\]
\be
\leq \eps_n\norm{(v_n\cdot\nb)v_n}_{L^{\kappa,\lambda}(\tilde{Q}^+)}
\leq N\eps_n.
\label{14feb100}
\ee
Next we define $\tilde{v}_n$ and $\tilde{q}_n$ as follows:
\[
\tilde{v}_n=v_n-\hat{v}_n,\qquad \tilde{q}_n=q_n-\hat{q}_n.
\]
Then it is straightforward that $\tilde{v}_n$ and $\tilde{q}_n$ solve
\[
\pr_s \tilde{v}_n-\Delta \tilde{v}_n+\nb \tilde{q}_n=0,\quad
{\rm{div}}\,\,\tilde{v}_n=0\quad \mbox{ in }\tilde{Q}^+,
\]
\[
\tilde{v}_n=0\quad\mbox{ on }\bke{\tilde{B}^+\cap\{x_3=0\}}
\times \bkt{-\bkp{\frac{3}{4}}^2,0}
\]
and  $\tilde{v}_n, \tilde{q}_n$ satisfy the following estimate
\[
\norm{\nb^2\tilde{v}_n}_{L{\kappa,\lambda}(Q^+_{\frac{5}{8}})}
+\norm{\nb\tilde{q}_n}_{L{\kappa,\lambda}(Q^+_{\frac{5}{8}})}
\leq N(1+\eps_n),
\]
and furthermore, 
for $\tilde{\kappa}$ with $3/\tilde{\kappa}+2/\lambda=1$ we obtain
\be
\norm{\nb^2\tilde{v}_n}_{L{\tilde{\kappa},\lambda}(Q^+_{\frac{9}{16}})}
+\norm{\nb\tilde{q}_n}_{L{\tilde{\kappa},\lambda}(Q^+_{\frac{9}{16}})}
\leq N(1+\eps_n ).
\label{14feb200}
\ee
Those estimates are again due to the local boundary 
estimate for the Stokes system (see \cite[Proposition 1-2]{S00}).
By the Poincar\'e inequality, we have
\[
\tilde{D}(q_n,\theta)\leq 
N\left(\tilde{D}_1(\hat{q}_n,\theta)+\tilde{D}_1(\tilde{q}_n,\theta)\right)
\]
We note that $\tilde{D}_1(\hat{q}_n,\theta)$ goes to zero 
as $n\,\raro\,\infty$ because of the estimate \eqref{14feb100}.
On the other hand, using the H\"older inequality and the estimate 
\eqref{14feb200}, we can show
\[
\tilde{D}_1(\tilde{q}_n,\theta)
=\frac{1}{\theta}\bke{\int_{-\theta^2}^0\bke{\int_{B^+(\theta)}
\abs{\nb\tilde{q}_n}^{\kappa}dy}
^{\frac{\lambda}{\kappa}}ds}^{\frac{1}{\lambda}}
\]
\[
\leq \theta^2\bke{\int_{-\theta^2}^0\bke{\int_{B^+(\theta)}
\abs{\nb\tilde{q}_n}^{\tilde{\kappa}}dy}
^{\frac{\lambda}{\tilde{\kappa}}}ds}^{\frac{1}{\lambda}}
\leq N \theta^2(1+\eps_n).
\]
So summing up, we obtain
\be
\liminf_{n\raro\infty}\tilde{D}(q_n,\theta)\leq 
\lim_{n\raro\infty}N_2\theta^2(1+\eps_n)
\leq N_2\theta^{1+\ap_0},
\label{contra:200}
\ee
where $N_2$ is an absolute constant.
At the beginning in \eqref{weakq:100} 
we can take an absolute constant $N$ bigger than $2(N_1+N_2)$ where
$N_1$ and $N_2$ are absolute constants in 
\eqref{hs:100} and \eqref{contra:200}, respectively.
Then this leads to a contradiction since
\[
2(N_1+N_2)\theta^{1+\ap_0}\leq N\theta^{1+\ap_0}\leq 
\liminf_{n\raro\infty}\psi_n(\theta)\leq (N_1+N_2)\theta^{1+\ap_0}. 
\]
This completes the proof.
\end{pf}
The Lemma above is the main part of the Lemma \ref{mod:lemma}.
Since the rest of the proof of Lemma \ref{mod:lemma} can be achieved by 
following similar procedures in \cite{S02},
we present only a brief sketch of the main idea of Lemma \ref{mod:lemma}.
\begin{lpf}
We first note that the Lemma above allows iterations
(compare \cite[Lemma 4.2]{S02}), and therefore,
there exists a positive constant $\ap_1 < 1$ 
\sidenote{removed $\gamma$ dependence in notation}
such that (compare \cite[Lemma 4.3 and Lemma 4.4]{S02})
\[
\left(C^{\frac{1}{3}}(r)+\tilde{D}(r)\right)
\leq N\bke{\frac{r}{\rho}}^{1+\ap_1}
\left(C^{\frac{1}{3}}(\rho)+\tilde{D}(\rho) 
+ m _\gamma(f) \rho^{\beta +1}
\right), 
\quad (r \le \rho).
\]
We consider for any $w\in\overline{Q^+_{z,r_1/2}}$
\[
\tilde{C}(w,r):=\frac{1}{r^2}\int_{Q_{w,r}\cap Q^+_{z,r_1/2}}
\abs{u-(u)_{a}}^3 dz,\quad
(u)_a=\aint_{Q_{w,r}\cap Q^+_{z,r_1/2}} u(z)dz,
\] 
and we can show that for any $r<r_1/4$
\[
\tilde{C}^{\frac{1}{3}}(w,r)\leq Nr^{1+\ap_1}, 
\]
where $N=N(\lambda,\gamma, \norm{f}_{M^{2,\gamma}})$ is an absolute constant. 
This argument can be proved using the same method as Lemma 5.2 in \cite{S02},
and therefore we omit the details. 
The regularity of $u$ at $z$ is a standard consequence of this estimate.
This completes the sketch of the proof.
\qed
\end{lpf}


\section{Local boundary regularity} 

In this section, we will present the proof of our main theorem
(see Theorem \ref{TH1}).
We first begin with an estimate for the scaled $L^3$-norm of 
suitable weak solutions.

\begin{lemm}\label{basiclemma}
Suppose (without loss of generality) $z = (0,0)$.  
Let $p, q$ be the numbers in \eqref{reg:20} and 
$Q^+_r=B^+_r\times (-r^2,0)$.
Suppose $u$ is a suitable weak solution
of the Navier-Stokes equations 
satisfying Definition \ref{defsws:100}. If
$u\in L^{p,q}(Q^+_r)$ and
$u=0$ on $(B_r\cap\{x_3=0\})\times (-r^2,0)$, then
\be
C(r)\leq NA^{\frac{1}{q}}(r)E^{1-\frac{1}{q}}(r)G(r).
\label{l3est:100}
\ee
\end{lemm}
\begin{pf}
We take $\ap, \beta,$ and $\delta$ such that 
$\ap=1/q, \beta=(1/3)(1-1/q)$, and $\delta=1/p$. 
Thus $2\ap+6\beta+p\delta=3$, and, therefore,
using the H\"older inequality, we obtain
\[
\int_{B^+_r}\abs{u}^3dx\leq \bke{\int_{B^+_r}\abs{u}^2dx}^{\ap}
\bke{\int_{B^+_r}\abs{u}^6dx}^{\beta}
\bke{\int_{B^+_r}\abs{u}^pdx}^{\delta}
\] 
\[
\leq N\bke{\int_{B^+_r}\abs{u}^2dx}^{\ap}
\bke{\int_{B^+_r}\abs{\nb u}^2dx}^{3\beta}
\bke{\int_{B^+_r}\abs{u}^pdx}^{\delta},
\]
where Sobolev embedding is used.
Integrating in time, we obtain
\[
\int_{Q^+_r}\abs{u}^3dz\leq 
N\bke{\sup_{-r^2\leq t\leq 0}\int_{B^+_r}\abs{u}^2dx}^{\ap}
\int_{-r^2}^0\bke{\int_{B^+_r}\abs{\nb u}^2dx}^{3\beta}
\bke{\int_{B^+_r}\abs{u}^pdx}^{\delta}dt
\]
\[
\leq 
N\bke{\sup_{-r^2\leq t\leq 0}\int_{B^+_r}\abs{u}^2dx}^{\ap}
\bke{\int_{Q^+_r}\abs{\nb u}^2dz}^{3\beta}
\bke{\int_{-r^2}^0\bkp{\int_{B^+_r}
\abs{u}^pdx}^{\frac{q}{p}}dt}^{\frac{1}{q}},
\]
where we used $3\beta+\delta=1$ and H\"older inequality.
Dividing both sides by $r^2$, we obtain \eqref{l3est:100}.
This completes the proof.
\end{pf}

\begin{rem}
The above estimate \eqref{l3est:100} is also true in the case
$1<q\leq\infty$ and $1<p<\infty$, although we 
restrict to numbers $p, q$ satisfying \eqref{reg:20}
\qed
\end{rem}

An immediate consequence of the local energy inequality is 
\be
A(\frac{r}{2})+E(\frac{r}{2})\leq N\left(C^{\frac{2}{3}}(r)+C(r)
+G(r)\tilde{D}(r)+r\int_{Q^+_r}|f|^2dz\right), 
\label{localenergy:200}
\ee
\[
\leq N\left(C^{\frac{2}{3}}(r)+C(r)+G(r)\tilde{D}(r)
+r^{2(\gamma+1)}m^2_{\gamma}\right).
\]
For those exponents $\kappa$ and $\lambda$ in 
\eqref{kalam:10} and \eqref{kalampq:100}, we can show (compare \eqref{preest:890})
\begin{equation}
\norm{(u\cdot\nb) u}_{L^{\kappa,\lambda}(Q^+_\rho)}
\leq N\rho E^{\frac{1}{\lambda}}(\rho)A^{\frac{3-2\kappa}{2\kappa}}(\rho).
\label{kalam:300}
\end{equation}
Its verification is straightforward, and we omit the details.

In next lemma we prove an estimate for the pressure.
\begin{lemm}\label{preest:100}
Suppose $z=(x,t), x\in\Gamma, t-\rho^2>0,$ 
and $t<T$. Then for $0 \leq r \leq \rho/4$,
\be
\tilde{D}_1(r)\leq N\left(\bkp{\frac{\rho}{r}}(E^{\frac{1}{\lambda}}(\rho)
           A^{\frac{3-2\kappa}{2\kappa}}(\rho)+\rho^{\gamma+1}m_{\gamma})
         +\bkp{\frac{r}{\rho}}(E^{\frac{1}{2}}(\rho)+\tilde{D}_1(\rho))\right),
\label{pres:100}
\ee
where $\kappa$ and $\lambda$ are numbers in 
\eqref{kalam:10} and \eqref{kalampq:100}.
\end{lemm}
\begin{pf}
Without loss of generality, we assume $x=0$.
We choose a domain $\tilde{B}^+$ with a smooth boundary such that
$B^+_{\rho/2}\subset\tilde{B}^+\subset B^+_{\rho}$, and we denote
$\tilde{Q}^+:=\tilde{B}^+\times (t-\rho^2,t)$.
We note first that, by the definition of $m_{\gamma}$ and the H\"older 
inequality, we have
\be
\norm{f}_{L^{\kappa,\lambda}(Q^+_\rho)}\leq N\rho^{\gamma+2}m_{\gamma},
\quad \norm{\nb u}_{L^{\kappa,\lambda}(Q^+_\rho)}
\leq N\rho^2E^{\frac{1}{2}}(\rho).
\label{morrey:1000}
\ee
Let $v$ and $p_1$ be the unique solution to the following initial boundary
value problem for the Stokes system
\[ 
\pr_t v-\Delta v+\nb p_1=(u \cdot \nb)u+f,\quad
\rm{div}\,\, v=0\qquad
\mbox{ in }\,\,\tilde{Q}^+,
\]
\[
(p_1)_{\tilde{B}^+}(t)=\aint_{\tilde{B}^+}p_1(y,t)\ind y=0,\quad
t\in (t-\rho^2,t)
\]
\[
v=0\quad\mbox{ on }\pr \tilde{B}^+\times [t-\rho^2,t],
\]
\[
v=0\quad\mbox{ on }\tilde{B}^+\times \{t=t-\rho^2\}.
\]
Then $v$ and $p_1$ satisfy the following estimate 
(see \cite[Theorem 3.1]{GS})
\[
\begin{split}
\frac{1}{\rho^2} & \norm{v}_{L^{\kappa,\lambda}(\tilde{Q}^+)}
+\frac{1}{\rho}\norm{\nb v}_{L^{\kappa,\lambda}(\tilde{Q}^+)}
+\norm{\pr_t v}_{L^{\kappa,\lambda}(\tilde{Q}^+)}
+\norm{\nb^2 v}_{L^{\kappa,\lambda}(\tilde{Q}^+)} \\
&+\frac{1}{\rho}\norm{p_1}_{L^{\kappa,\lambda}(\tilde{Q}^+)}
+\norm{\nb p_1}_{L^{\kappa,\lambda}(\tilde{Q}^+)} \\
&\leq N\bke{\norm{u\nb u}_{L^{\kappa,\lambda}(\tilde{Q}^+)}
+\norm{f}_{L^{\kappa,\lambda}(\tilde{Q}^+)}} 
\leq N\bke{\norm{u\nb u}_{L^{\kappa,\lambda}(Q^+_{\rho})}
+\norm{f}_{L^{\kappa,\lambda}(Q^+_{\rho})}} \\
&\leq N\bke{\rho E^{\frac{1}{\lambda}}(\rho)A^{\frac{3-2\kappa}{2\kappa}}(\rho)
+\rho^{\gamma+2}m_{\gamma}},
\end{split}
\]
where we used \eqref{kalam:300} and \eqref{morrey:1000}.

Let $w=u-v$ and $p_2=p-(p)_{B^+_{\rho/2}}-p_1$. Then $w,p_2$
solve the following boundary value problem:
\[
\pr_t w-\Delta w+\nb p_2=0,\quad
\rm{div}\,\, w=0
\qquad\mbox{ in }\,\,\tilde{Q}^+,
\]
\[
w=0\quad\mbox{ on }\,\,\bke{\pr\tilde{B}^+\cap\{x_3=0\}}
\times \bkt{t-\rho^2,t}
\]
Now we take $\kappa'$ ($\kappa'$ is different than $\kappa^*$) such that 
$3/\kappa'+2/\lambda=2$.
Then from the local estimate near the boundary for the 
Stokes system (see \cite[Proposition 2]{S00}), we obtain
\[
\norm{\nb^2 w}_{L^{\kappa',\lambda}(Q^+_{\frac{\rho}{4}})}
+\norm{\nb p_2}_{L^{\kappa',\lambda}(Q^+_{\frac{\rho}{4}})}
\]
\[
\leq \frac{N}{\rho^2}\bke{\frac{1}{\rho^2}\norm{w}_{L^{\kappa,\lambda}
(Q^+_{\frac{\rho}{2}})}
+\frac{1}{\rho}\norm{\nb w}_{L^{\kappa,\lambda}(Q^+_{\frac{\rho}{2}})} 
+\frac{1}{\rho}\norm{p_2}_{L^{\kappa,\lambda}(Q^+_{\frac{\rho}{2}})}}
\equiv \frac{N}{\rho^2}I 
\]
Using Sobolev imbedding, the right side can be estimated as follows:
\[
I\leq \left(\frac{1}{\rho}\norm{\nb u}_{L^{\kappa,\lambda}
(Q^+_{\frac{\rho}{2}})} 
+\norm{\nb p}_{L^{\kappa,\lambda}(Q^+_{\frac{\rho}{2}})}
+\frac{1}{\rho}\norm{\nb v}_{L^{\kappa,\lambda}(Q^+_{\frac{\rho}{2}})} 
+\frac{1}{\rho}\norm{p_1}_{L^{\kappa,\lambda}(Q^+_{\frac{\rho}{2}})}\right).
\]
Due to the second inequality in \eqref{morrey:1000}, we obtain
\[
\norm{\nb p_2}_{L^{\kappa',\lambda}(Q^+_{\frac{\rho}{4}})}\leq 
\frac{N}{\rho^2}\left(\rho E^{\frac{1}{2}}(\rho)+\rho \tilde{D}_1(\rho)
+\rho E^{\frac{1}{\lambda}}(\rho)A^{\frac{3-2\kappa}{2\kappa}}(\rho)
+\rho^{\gamma+2}m_{\gamma}\right)
\]
\[
=\frac{N}{\rho}\left(E^{\frac{1}{2}}(\rho)+\tilde{D}_1(\rho)
+E^{\frac{1}{\lambda}}(\rho)A^{\frac{3-2\kappa}{2\kappa}}(\rho)
+\rho^{\gamma+1}m_{\gamma}\right).
\]
Now we assume $0\leq r\leq \rho/4$.
Noting that $||\nb p_2||_{L^{\kappa,\lambda}(Q^+_r)}
\leq Nr^2||\nb p_2||_{L^{\kappa',\lambda}(Q^+_r)}$, we have
\[
\tilde{D}_1(r)=\frac{1}{r}\norm{\nb p}_{L^{\kappa,\lambda}(Q^+_r)}
\leq \frac{1}{r}\bke{\norm{\nb p_1}_{L^{\kappa,\lambda}(Q^+_r)}
                 +\norm{\nb p_2}_{L^{\kappa,\lambda}(Q^+_r)}}
\]
\[
\leq \frac{1}{r}\bke{\norm{\nb p_1}_{L^{\kappa,\lambda}(Q^+_{\rho})}
                 +r^2\norm{\nb p_2}_{L^{\kappa',\lambda}(Q^+_r)}}
\]
\[
\leq N(\frac{\rho}{r})\left(E^{\frac{1}{\lambda}}(\rho)
A^{\frac{3-2\kappa}{2\kappa}}(\rho)+\rho^{\gamma+1}m_{\gamma}\right)+
\]
\[
+N(\frac{r}{\rho})\left(E^{\frac{1}{2}}(\rho)+\tilde{D}_1(\rho)
+E^{\frac{1}{\lambda}}(\rho)A^{\frac{3-2\kappa}{2\kappa}}(\rho)
+\rho^{\gamma+1}m_{\gamma}\right)
\]
\[
\leq N(\frac{\rho}{r})\left(E^{\frac{1}{\lambda}}(\rho)
A^{\frac{3-2\kappa}{2\kappa}}(\rho)+\rho^{\gamma+1}m_{\gamma}\right)
+N(\frac{r}{\rho})\left(E^{\frac{1}{2}}(\rho)+\tilde{D}_1(\rho)\right).
\]
This completes the proof.
\end{pf}

Now we are ready to present the proof of Theorem \ref{TH1}.
\begin{sol}\,\,\,
We recall first, due to Lemma \ref{basiclemma},
\be
C(r)\leq NA^{\ap}(r)E^{3\beta}(r)G(r),\quad \ap=\frac{1}{q},\,\,
\beta=\frac{q-1}{3q},
\label{main:1000}
\ee
and, due to Sobolev imbedding, we have 
\be
\tilde{D}(r)\leq N_d\tilde{D}_1(r)
\label{const:nd}
\ee

Let $4r<\rho$. We consider $C(r)+\tilde{D}_1(r)$.
Recalling the estimate \eqref{pres:100} for the pressure, we obtain 
\[
\tilde{D}_1(r)+C(r)\leq 
NA^{\ap}(r)E^{3\beta}(r)G(r)+
\]
\[
+N\bkp{\frac{\rho}{r}}\bke{E^{\frac{1}{\lambda}}(\frac{\rho}{4})
A^{\frac{3-2\kappa}{2\kappa}}(\frac{\rho}{4})+m_{\gamma}\rho^{\gamma+1}}
+N\bkp{\frac{r}{\rho}}\bke{E^{\frac{1}{2}}(\frac{\rho}{4})
+\tilde{D}_1(\frac{\rho}{4})}
\]
\[
\equiv I+II+III.
\]
We first consider the first term $I$.
Since $\ap+3\beta=1$, by using the local energy inequality 
\eqref{localenergy:200}, we have
\[
I\leq N\left(C^{\frac{2}{3}}(2r)+C(2r)+G(2r)\tilde{D}(2r)
+r^{2(\gamma+1)}m^2_{\gamma}\right)G(r)
\]
\[
\leq N\bke{\bkp{\frac{\rho}{r}}^{\frac{7}{3}}C^{\frac{2}{3}}(\rho)G(\rho)
+\bkp{\frac{\rho}{r}}^3C(\rho)G(\rho)
+\bkp{\frac{\rho}{r}}^3G^2(\rho)\tilde{D}_1(\rho)
+\bkp{\frac{\rho}{r}}r^{2(\gamma+1)}m^2_{\gamma}G(\rho)}
\]
\be
\leq N\bke{\bkp{\frac{\rho}{r}}^3C(\rho)G(\rho)
+\bkp{\frac{\rho}{r}}^3G^2(\rho)\tilde{D}_1(\rho)
+\bkp{\frac{\rho}{r}}G(\rho)
+\bkp{\frac{\rho}{r}}\rho^{2(\gamma+1)}m^2_{\gamma}G(\rho)},
\label{main:150}
\ee
where we used Young's inequality and 
\[
C(2r)\leq N\bkp{\frac{\rho}{r}}^2C(\rho),\quad
\tilde{D}_1(2r)\leq N\bkp{\frac{\rho}{r}}\tilde{D}_1(\rho),\quad
G(2r)\leq N\bkp{\frac{\rho}{r}}G(\rho).
\]
For the third term $III$, 
again using energy inequality \eqref{localenergy:200}, 
we have
\[
III
\leq N\bkp{\frac{r}{\rho}}\bke{\bke{C^{\frac{1}{3}}(\frac{\rho}{2})
+C^{\frac{1}{2}}(\frac{\rho}{2})
+G^{\frac{1}{2}}(\frac{\rho}{2})
\tilde{D}_1^{\frac{1}{2}}(\frac{\rho}{2})+m_{\gamma}\rho^{\gamma+1}}
+\tilde{D}_1(\frac{\rho}{2})}
\]
\[
\leq N\bkp{\frac{r}{\rho}}
\left(G(\frac{\rho}{2})+C^{\frac{1}{3}}(\frac{\rho}{2})
+C(\frac{\rho}{2})+\tilde{D}_1(\frac{\rho}{2})+m_{\gamma}\rho^{\gamma+1}
\right), 
\]
where we used Young's inequality, i.e. $ab\leq a^l/l+b^m/m$ where 
$1/l+1/m=1, 1<l,m<\infty$. By~\eqref{main:150}, note that 
\[
C^{\frac{1}{3}}(\frac{\rho}{2})\leq 
NA^{\frac{\ap}{3}}(\frac{\rho}{2})E^{\beta}(\frac{\rho}{2})
G^{\frac{1}{3}}(\frac{\rho}{2})
\]
\[
\leq N\left(C^{\frac{2}{3}}(\rho)+C(\rho)
+G(\rho)\tilde{D}_1(\rho)+m^2_{\gamma}\rho^{2(\gamma+1)}
\right)^{\frac{1}{3}}G^{\frac{1}{3}}(\rho)
\]
\[
\leq N\left(C^{\frac{2}{9}}(\rho)+C^{\frac{1}{3}}(\rho)
+G^{\frac{1}{3}}(\rho)\tilde{D}^{\frac{1}{3}}_1(\rho)
+m^{\frac{2}{3}}_{\gamma}\rho^{\frac{2}{3}(\gamma+1)}
\right)G^{\frac{1}{3}}(\rho).
\]
Again applying Young's inequality, 
we obtain
\[
C^{\frac{1}{3}}(\frac{\rho}{2})
\leq N\left(C(\rho)+\tilde{D}_1(\rho)+
G^{\frac{3}{7}}(\rho)+G(\rho)+
m^{\frac{2}{3}}_{\gamma}\rho^{\frac{2}{3}(\gamma+1)}
G^{\frac{1}{3}}(\rho)
\right).
\]
Summing up, we obtain
\be
III\leq N(\frac{r}{\rho})\bke{
C(\rho)+\tilde{D}_1(\rho)+G(\rho)
+G^{\frac{3}{7}}(\rho)
+m_{\gamma}\rho^{\gamma+1}
+m^{\frac{2}{3}}_{\gamma}\rho^{\frac{2}{3}(\gamma+1)}
G^{\frac{1}{3}}(\rho)}.
\label{main:250}
\ee
It remains to consider the second term $II$.
Since $1/\lambda+(3-2\kappa)/(2\kappa)=1$, by \eqref{localenergy:200},
we obtain 
\[
II\leq N(\frac{\rho}{r})\left( C^{\frac{2}{3}}(\frac{\rho}{2})
+C(\frac{\rho}{2})+G(\frac{\rho}{2})\tilde{D}_1(\frac{\rho}{2})
+\rho^{2(\gamma+1)}m^2_{\gamma}+m_{\gamma}\rho^{\gamma+1}\right).
\]
Using the same procedure as above, using~\eqref{main:1000}
and Young's inequality, we obtain 
\[
C^{\frac{2}{3}}(\frac{\rho}{2})
\leq N\left(C^{\frac{4}{9}}(\rho)+C^{\frac{2}{3}}(\rho)
+G^{\frac{2}{3}}(\rho)\tilde{D}^{\frac{2}{3}}_1(\rho)
+m^{\frac{4}{3}}_{\gamma}\rho^{\frac{4}{3}(\gamma+1)}
\right)G^{\frac{2}{3}}(\rho)
\]
\[
\leq N\Big(\bkp{G(\rho)+G^{\frac{1}{2}}(\rho)}C(\rho)
+G(\rho)\tilde{D}_1(\rho)+
\]
\[
+G^2(\rho)+G^{\frac{4}{5}}(\rho)
+m^{\frac{4}{3}}_{\gamma}\rho^{\frac{4}{3}(\gamma+1)}
G^{\frac{2}{3}}(\rho)\Big),
\]
and 
\[
C(\frac{\rho}{2})
\leq N\left(C^{\frac{2}{3}}(\rho)+C(\rho)
+G(\rho)\tilde{D}_1(\rho)+m^2_{\gamma}\rho^{2(\gamma+1)}
\right)G(\rho)
\]
\[
\leq N\left(C(\rho)G(\rho)+G^2(\rho)\tilde{D}_1(\rho)
+G(\rho)+m^2_{\gamma}\rho^{2(\gamma+1)}
G(\rho)\right).
\]
Summing up all together, we have
\[
II\leq N(\frac{\rho}{r})\Big(
\bkp{G(\rho)+G^{\frac{1}{2}}(\rho)}C(\rho)
+\bkp{G(\rho)+G^2(\rho)}\tilde{D}_1(\rho)+
\]
\be
G^{\frac{4}{5}}(\rho)+G^2(\rho)
+m^{\frac{4}{3}}_{\gamma}\rho^{\frac{4}{3}(\gamma+1)}
G^{\frac{2}{3}}(\rho)
+\rho^{2(\gamma+1)}m^2_{\gamma}+m_{\gamma}\rho^{\gamma+1}\Big).
\label{main:350}
\ee
Adding \eqref{main:150}, \eqref{main:250}, and \eqref{main:350}, we obtain
\[
C(r)+\tilde{D}_1(r)\leq 
N\bke{(\frac{\rho}{r})^3G(\rho)+(\frac{\rho}{r})
\bkp{G(\rho)+G^{\frac{1}{2}}(\rho)}+(\frac{r}{\rho})}C(\rho)+
\]
\[
+N\bke{(\frac{\rho}{r})^3G(\rho)+(\frac{\rho}{r})
\bkp{G(\rho)+G^2(\rho)}+(\frac{r}{\rho})}\tilde{D}_1(\rho)+
\]
\[
+(\frac{\rho}{r})(G^{\frac{3}{7}}(\rho)+G^2(\rho))+
(\frac{\rho}{r})m_{\gamma}\rho^{\gamma+1}
+(\frac{\rho}{r})m^2_{\gamma}\rho^{2(\gamma+1)}G(r)+
\]
\[
+(\frac{\rho}{r})m^2_{\gamma}\rho^{2(\gamma+1)}
+m^{\frac{2}{3}}_{\gamma}\rho^{\frac{2}{3}(\gamma+1)}
G^{\frac{1}{3}}(\rho)+m^{\frac{4}{3}}_{\gamma}\rho^{\frac{4}{3}(\gamma+1)}
G^{\frac{2}{3}}(\rho).
\]
We first choose $\theta\in [0,1/2]$ such that $N\theta<1/4$ where 
$N$ is an absolute constant in the above inequality.
By replacing $r,\rho$ by $\theta r$ and $r$, we obtain
\[
C(\theta r)+\tilde{D}_1(\theta r)
\leq N\Big(\bkp{\frac{1}{\theta^3}G(r)+\frac{1}{\theta}
G^{\frac{1}{2}}(r)+\theta}C(r)+
\]
\be
+\bkp{\frac{1}{\theta^3}G(r)+\frac{1}{\theta}
G^2(r)+\theta}\tilde{D}_1(r)+\phi(r)\Big),
\label{main:2000}
\ee
where
\[
\phi(r)=\frac{1}{\theta}\bke{G^{\frac{3}{7}}(r)+G^2(r)}+
\frac{1}{\theta}m_{\gamma}r^{\gamma+1}
+\frac{1}{\theta}m^2_{\gamma}r^{2(\gamma+1)}G(r)+
\]
\[
+\frac{1}{\theta}m^2_{\gamma}r^{2(\gamma+1)}
+m^{\frac{2}{3}}_{\gamma}r^{\frac{2}{3}(\gamma+1)}
G^{\frac{1}{3}}(r)+m^{\frac{4}{3}}_{\gamma}r^{\frac{4}{3}(\gamma+1)}
G^{\frac{2}{3}}(r).
\]
Now we fix $r_0<\min\{1,\theta\eps^3/(1+m_{\gamma})\}$ 
such that for all $r\leq r_0$ 
\[
G(r)<\min\bket{\frac{\theta^3}{2^{10}N}, \frac{\theta^2}{2^{10}N},
\frac{\eps^3\theta}{2^{10}N(N_d+1)m^2_{\gamma}}, 
(\frac{\eps^3\theta}{2^{10}N(N_d+1)})^{\frac{7}{3}},
\frac{\eps^9}{(2^{10}N(N_d+1))^3m^2_{\gamma}}},
\]
where $N, N_d$ are absolute constants in \eqref{main:2000} 
and \eqref{const:nd}, respectively, and 
$\eps$ is the fixed positive number in Lemma \ref{mod:lemma}.
Then one can check that $\phi(r)<\eps^3/64N(N_d+1)$
and moreover, we can show that for any $r<r_0$
\[
C(\theta r)+\tilde{D}_1(\theta r)\leq \frac{1}{2}\bke{C(r)+\tilde{D}_1(r)}
+\phi(r)
\]
By iterating, we have
\[
C(\theta^k r)+\tilde{D}_1(\theta^k r)\leq (\frac{1}{2})^k 
\bke{C(r)+\tilde{D}_1(r)}
+\sum_{i=0}^{k-1}\frac{N}{2^{k-1-i}}\phi(\theta^i r).
\]
\be
\leq (\frac{1}{2})^k \bke{C(r)+\tilde{D}_1(r)}+\frac{\eps^3}{64(N_d+1)}
\label{it:10}
\ee
If $z$ is a singular point, then there exists $r_1>0$
such that $C^{\frac{1}{3}}(r)+D(r)\geq \eps$ for every $r\leq r_1$ by 
Lemma \ref{mod:lemma}.
However, 
this leads to a contradiction since for a sufficiently small $r_2<r_1$
\[
C(r_2)+\tilde{D}(r_2)\leq C(r_2)+N_d\tilde{D}_1(r_2)
\leq\frac{\eps^3}{64},
\]
which immediately implies that 
$C^{\frac{1}{3}}(r)+D(r)\leq \eps/2$.
This completes the proof.
\qed
\end{sol}

The following is a direct consequence of Theorem \ref{TH1}.
\begin{cor}\label{swecor}
Let $u$ be a weak solution of the Navier-Stokes equations 
satisfying Definition \ref{defsws:100}.
Assume further that $z=(x,t)\in\Gamma\times I$ and 
for some $r_0>0$ 
\begin{equation}
u\in L^{r,s}(Q^+_{z,r_0}),
\qquad
\frac{3}{r}+\frac{2}{s}=1,\,\,3<r<\infty.
\label{SPC:sec3}
\end{equation}
Then $z$ is a regular point.
\end{cor}
\begin{pf}
We observe that, as mentioned in the introduction, 
$u\in L^4(Q^+_{z,r_0})$ for weak solutions satisfying \eqref{SPC:sec3}.
To be more precise, we can show by interpolation that 
\[
\norm{u}_{L^4(Q^+_{z,r_0})}\leq \norm{u}^{\frac{1}{2}}_{L^{r,s}(Q^+_{z,r_0})}
\norm{u}^{\frac{1}{q}}_{L^{2,\infty}(Q^+_{z,r_0})}
\norm{u}^{\frac{3}{2p}}_{L^{6,2}(Q^+_{z,r_0})}.
\]
The above estimate is true even in case $(r,s)=(3,\infty)$ or 
$(r,s)=(\infty,2)$, although our analysis does not include 
such cases.
We conclude by the above estimate that 
$u$ is a suitable weak solution in $Q^+_{z,r_0}$, namely
$u$ satisfies the local energy inequality \eqref{lei} in a neighborhood 
of $z$.
We also note that there exists a number $\tilde{r}$ 
such that $3/\tilde{r}+2/s=2$ and, by the H\"older inequality, we have 
\[
\frac{1}{\rho}\norm{u}_{L^{\tilde{r},s}(Q^+_{z,\rho})}
\leq C\norm{u}_{L^{r,s}(Q^+_{z,\rho})},\qquad \mbox{ for any }\rho\leq r_0.
\]
Since the right hand side above is finite, and it can be arbitrary small
for sufficiently small $\rho$ 
by assumption \eqref{SPC:sec3}, the condition \eqref{nse:abs30}
in Theorem \ref{TH1} is satisfied. 
This completes the proof.
\end{pf}

\begin{rem}
The condition \eqref{SPC:sec3}, including the case $(r,s)=(\infty,2)$,
is called a \ps condition. 
It is known that weak solutions of the Navier-Stokes equations
are locally regular at an interior point provided that a \ps condition
is assumed near the point (see \cite{S62} and \cite{S88}).
This result, recently, was extended up to the boundary 
(see \cite{k_bnse} and \cite{VS02}), and, therefore,
Corollary \ref{swecor} is already implied by \cite{k_bnse} and \cite{VS02}. 
Our regularity criterion \eqref{nse:abs30}, 
however, gives a simple proof of the regularity of weak solutions
near the boundary under a \ps condition (although the case $(r,s)=(\infty,2)$
is not covered by our analysis).
\end{rem}


\section{Partial regularity}

In this section, as an application of Theorem \ref{TH1},
we investigate the size of the possible singular set 
under additional integrability assumptions on weak solutions
(see Assumption \ref{assume100} below).
As we saw in Corollary \ref{swecor},
we have a simple proof for weak solutions that  
a \ps condition implies regularity up to the boundary. 
It is, however, an open question whether
or not weak solutions (or suitable weak solutions) satisfy
the \ps conditions.
It was proved that the size of a possible 
singular set for suitable weak solutions 
is of $1$-dimensional Hausdorff measure zero (see 
\cite{CKN} and \cite{S02} for the interior case 
and for the boundary case, respectively).
We remark that in \cite{CL00} the estimate of the Hausdorff 
measure of the singular set for suitable weak solutions
was improved by a logarithmic factor
for the interior case.
 
Our aim in this section is to present the proof of Theorem \ref{TH2}, 
which says that the size of singular set for weak solutions
can be reduced under additional integrability assumptions, which are  
weaker than \ps conditions.
We note, however, that our result is weaker than what one 
gets from the \ps conditions; in that case full regularity is implied, 
but in our case we have partial regularity. 

We start with recalling the following condition 
that is assumed in Theorem \ref{TH2} for weak solutions.

\begin{assume}\label{assume100}
Let $u$ be a weak solution of the Navier-Stokes equations 
satisfying Definition \ref{defsws:100} and 
\begin{equation*}
u\in L^{m,l}(Q)=L^m\bke{(0,T);L^l(\om)},
\end{equation*}
where either 
\begin{equation}
1<\frac{3}{l}+\frac{2}{m},\quad \frac{2}{l}+\frac{2}{m}<1, \quad
l> m,
\label{ass200;100}
\end{equation}
or 
\begin{equation}
1<\frac{3}{l}+\frac{2}{m},\quad \frac{3}{l}+\frac{1}{m}<1, \quad
l\leq m.
\label{ass200;200}
\end{equation}
\qed
\end{assume}
We remark that although Assumption \ref{assume100}
is not justified by the formulation of weak solutions, 
it seems to be of independent interest to characterize the size 
of the singular set depending on the mixed norm $L^{l,m}$ of $u$.

We observe first that solutions satisfying Assumption \ref{assume100}
are in fact suitable weak solutions. 
This can be done by the interpolation argument
of Corollary \ref{swecor}. More precisely, we can show  
\[
\norm{u}_{L^4(Q)}\leq \norm{u}^{\sigma}_{L^{l,m}(Q)}
\norm{u}^{\frac{(6-\ap)}{2\ap}(1-\sigma)}_{L^{2,\infty}(Q)}
\norm{u}^{\frac{2}{\beta}(1-\sigma)}_{L^{6,2}(Q)},
\]
where
\[
\ap=\frac{\frac{3}{l}+\frac{2}{m}-\frac{5}{4}}{\frac{1}{l}+\frac{1}{2m}-\frac{3}{8}},\quad
\beta=\frac{4(\frac{3}{l}+\frac{2}{m}-\frac{5}{4})}{3(\frac{1}{l}+\frac{1}{m}-\frac{1}{2})},
\quad \sigma=\frac{l(4-l)}{4(l-\ap)}.
\]
Now we are ready to prove Theorem \ref{TH2}.

\begin{sol1}
We consider first the case that $l,m$ satisfy \eqref{ass200;200}.
For convenience, we define $\varsigma$ as
$\varsigma:=3/l+2/m-1$. 
For given $R>0$ we set $S_R:=S\cap B_R(0)$, where $S$ is the singular 
set of $u$. 
Let $z=(x,t)\in S_R$  and 
$\tilde{Q}_r(z)=\tilde{Q}_r(x,t):=Q_r(x,t)\cap Q$.
Using Theorem \ref{TH1} 
(boundary case) and Theorem \ref{int:mainthm400} (interior case)
in the Appendix, we see that there exists $\e_0$ such that 
\begin{equation}
\limsup_{r \to 0_+} r^{-\varsigma} 
\norm{u}_{L^m_tL^l_x\bke{\tilde{Q}_r(z)}} \ge \e_0.
\label{sin-ppt}
\end{equation}
It is clear that the Lebesgue measure of $S_R$ is zero, and so
we can choose an open 
bounded set $V\subset\R^3$ containing $\bar{S_R}$ with 
the volume of $V$ as small as we like. 
Moreover, due to \eqref{sin-ppt}, for any given $\delta > 0$,
and for any $z=(x,t)\in S_R$,
\sidenote{fixed this to match definition of Hausdorff measure}
there exists $r_z \in (0, \delta)$ such that 
$\norm{u}_{L^m_tL^l_x(\tilde{Q}_{r_z}(z))} \ge r^{\varsigma}_z\e_0$.
We denote by $\calO=\{\tilde{Q}_{r_z}(z):z\in S_R\}$ the collection of 
such open neighbourhoods.
We note that $\calO$ is an open covering of $S_R$, and, therefore, 
by a covering lemma (e.g. see \cite[Lemma 6.1]{CKN}), we can find 
a countable subfamily of disjoint cylinders 
$\tilde{Q}_j= \tilde{Q}_{r_{z_j}}(z_j), j\in J$ such that 
$S_R \subset \cup_{z \in S_R} \tilde{Q}_{r_z}(z)
\subset\cup_{j\in J} \tilde{Q}_{5r_j}(z_j)$.
Denoting $u_j(z)=u(z)$ if $z\in Q_j^+$, $u_j(z)=0$ otherwise,
we have
\[
\e_0^{m}\sum_{j\in J} r_j^{\varsigma m}\leq \sum_{j\in J} \int_I 
\bkp{\int_{\om} |u_j|^l
\ind x}^{m/l}\ind t 
\leq \int_I \Big(\sum_{j\in J}  \int_{\om} |u_j|^l
\ind x\Big)^{m/l} \ind t \leq  \norm{u}_{L^m_tL^l_x(V)}^m,
\]
where we used $l\leq m$ in the second inequality.
Since the volume of $V$ can be taken arbitrarily small,
as can $\delta$, we conclude that 
the $\varsigma m$-dimensional Hausdorff measure 
of $S_R$ is zero: ${\cal P} ^{\varsigma m} (S_R)=0$.
Since $R$ is arbitrary, we conclude that the singular set $S$
is of $\varsigma m=2-m+3m/l$ dimensional parabolic Hausdorff measure zero.
Observe that $0<\varsigma m<1$.

Next we consider the case $l,m$ satisfy \eqref{ass200;100}.
In this case we show, by interpolation, that 
\[
\norm{u}_{L^k(Q)}\leq \norm{u}^{1-\sigma}_{L^{l,m}(Q)}
\norm{u}^{\sigma}_{L^{2,\infty}(Q)},
\]
where 
\[
k=2(1-\frac{m}{l})+l\frac{m}{l}=2+m-\frac{2m}{l},\quad \sigma=\frac{m}{k}.
\]
It is clear that $4<k<5$ and $m/k\leq 1$, and 
therefore, due to the analysis of the case $l\leq m$, we conclude that
the singular set is at most of $5-k$ dimensional parabolic 
Hausdorff measure zero up to the boundary.
This completes the proof.
\qed
\end{sol1}

\section{Appendix}

In this Appendix we show that 
the regularity criterion \eqref{nse:abs30} holds also for the interior case,
whose proof requires slightly different estimates.
Since its verification for the interior case 
can be done by following a procedure similar to that of the boundary case,
with no significant difficulty,
we just present a sketch of the proof.
From now on we replace $Q^+_{z,r}$ and $B^+_{x,r}$ by 
$Q_{z,r}$ and $B_{x,r}$ in the scaling invariant functionals below, 
because we are concerned with local regularity at an interior point.
The interior case is in fact simpler than the boundary case, because 
the pressure is much easier to handle. 
We begin with the following lemma, which is analogous to 
Lemma \ref{lemma:seregin} of the boundary case. 

\begin{lemm}\label{int:lemma100}
Let $0<\theta<1/2$ and $0 < \beta < \gamma \le 2$.
There exist $\eps_2, r_2>0$ depending on  $\lambda,\theta, \gamma$ and $ \beta$
such that if 
$u$ is a suitable weak solution of the Navier-Stokes equations 
satisfying Definition \ref{defsws:100}, 
$z=(x,t)\in Q=\om\times I$ is an interior point, and 
$C^{\frac{1}{3}}(r)+\tilde{D}(r)+ m_{\gamma}(f) r^{\beta+1}<\eps_2$ 
for some $r\in (0,r_2)$, then
\[
\bke{C^{\frac{1}{3}}(\theta r)+\tilde{D}(\theta r)}
<N\theta^{1+\ap}\bke{C^{\frac{1}{3}}(r)+\tilde{D}(r)
+ m_{\gamma}(f) r^{\beta +1}},
\]
where $0<\ap<1$ and $N>0$ are absolute constants.
\end{lemm}
\begin{ispf}
Again assume $f=0$ for simplicity.
Suppose that the assertion is not true. 
Then there exist $z_n=(x_n,t_n)$, $r_n\searrow 0$,
and $\eps_n\searrow 0$ such that
$\varphi(r_n)=\eps_n$ but $\varphi(\theta r_n)>N\theta^{1+\ap}\eps_n$,
where $\varphi(r)=C^{\frac{1}{3}}(r)+\tilde{D}(r)$.
Using the change of variables $y=r^{-1}_n(x-x_n)$ and $s=r^{-2}_n(t-t_n)$,
we set $v_n(w):=\eps^{-1}_n r_nu(z)$ and 
$q_n(w):=\eps^{-1}_n r^2_n\bkp{p(z)-(p)_{B^+_{r_n}}(t)}$.
By the ``blow up'' procedure and compactness arguments, 
the limit equations become the Stokes system.
Since the pressure of the Stokes system system is harmonic in 
the spatial variables for the interior case, 
our arguments are much simpler than in the boundary case.
The other parts of the arguments are the same as 
in Lemma \ref{lemma:seregin}, and we omit the details.
\qed
\end{ispf}

Due to the Lemma \ref{int:lemma100}, we have 
the following lemma (compare to Lemma \ref{mod:lemma}).
Since the arguments are straightforward, we state it without proof. 

\begin{lemm}\label{int:lemma200}
There exists a constant $\eps>0$ depending on $\lambda,\gamma$ and
$\norm{f}_{M_{2,\gamma}}$ such that if $u$ is a suitable weak solution
of the Navier-Stokes equations satisfying Definition \ref{defsws:100},
$z=(x,t)\in Q=\om\times I$ is an interior point, and
\[
\liminf_{r\raro 0_+}\bke{C^{\frac{1}{3}}(r)+\tilde{D}(r)}<\eps,
\] 
then $z$ is a regular point.
\end{lemm}

Next we need the estimate for the pressure.
To do that, we observe that 
at an interior point, 
instead of \eqref{l3est:100} for the boundary case,
we can show 
\begin{equation}
C(r)\leq N\bke{A^{\frac{1}{q}}(r)E^{1-\frac{1}{q}}(r)G(r)
+A^{\frac{1}{2}}(r)G^2(r)}.
\label{appendix:l3}
\end{equation}
Indeed, 
\[
\int_{B_r}\abs{u}^3\ind y\leq N\left(\int_{B_r}\abs{u-(u)_a}^3\ind y
+\int_{B_r}\abs{(u)_a}^3\ind y\right),
\]
where $(u)_a=\aint_{B_r}u(y)\ind y$. We note that 
the first of the above
inequalities can be estimated as the same way as Lemma \ref{basiclemma}, 
and thus it is enough to consider the second one.
We observe that 
\[
\abs{(u)_a}\leq \frac{N}{r^{\frac{3}{2}}}
\bke{\int_{B_r}\abs{u(y)}^2\ind y}^{\frac{1}{2}}
\leq \frac{N}{r}A^{\frac{1}{2}}(r),\quad 
\abs{(u)_a}\leq\frac{N}{r^{\frac{3}{p}}}
\bke{\int_{B_r}\abs{u(y)}^p\ind y}^{\frac{1}{p}}.
\]
The second one is estimated as follows:
\[
\int_{B_r}\abs{(u)_a}^3\ind y\leq 
\frac{N}{r^{\frac{6}{p}-2}}A^{\frac{1}{2}}(r)
\bke{\int_{B_r}\abs{u(y)}^p\ind y}^{\frac{2}{p}}.
\]
Integrating in time, and using $q>2$, 
we obtain \eqref{appendix:l3}. Since the computations are straightforward,
the details are skipped.
In a similar manner, we also have the following estimate
in the interior case (compare to \eqref{kalam:300} for the boundary case).
\begin{equation}
\norm{(u\cdot\nb) u}_{L^{\kappa,\lambda}(Q_\rho)}
\leq N\rho \left(E^{\frac{1}{\lambda}}(\rho)
A^{\frac{3-2\kappa}{2\kappa}}(\rho)
+E^{\frac{1}{2}}(\rho)A^{\frac{2-\kappa}{2\kappa}}(\rho)
G^{\frac{\kappa-1}{\kappa}}(\rho)\right).
\label{appendix:unbu}
\end{equation}
Since its verification is similar to \eqref{appendix:l3}, we omit the details.

Using the above estimate, we have 
the pressure estimate, equivalent to Lemma \ref{preest:100}
of the boundary case.
Since the estimates \eqref{appendix:unbu} and 
\eqref{appendix:l3} are slightly different than  
\eqref{l3est:100} and \eqref{kalam:300} of the boundary case, 
the estimate of the pressure is slightly modified 
in the interior case. But it can be derived 
in the same manner, and so we skip its proof and just state it.
\begin{lemm}\label{int:lemma300}
Suppose that $z=(x,t)\in \om\times I$ is an interior point,
and $t-\rho^2>0,t<T$. Then for $0\leq r\leq \rho/4$,
\[
\tilde{D}_1(r)\leq N(\frac{\rho}{r})\left(E^{\frac{1}{\lambda}}(\rho)
           A^{\frac{3-2\kappa}{2\kappa}}(\rho)
+E^{\frac{1}{2}}(\rho)A^{\frac{2-\kappa}{2\kappa}}(\rho)
G^{\frac{\kappa-1}{\kappa}}(\rho)+\rho^{\gamma+1}m_{\gamma}\right)+
\]
\[
+N(\frac{r}{\rho})\left(E^{\frac{1}{2}}(\rho)+\tilde{D}_1(\rho)\right),
\]
where $\kappa$ and $\lambda$ are numbers satisfying \eqref{kalam:10} and
\eqref{kalampq:100}.
\end{lemm}

Using the estimate of the pressure in Lemma \ref{int:lemma300}, the
same regularity criterion for interior points can be proved as in
the boundary case. We have

\begin{thm}\label{int:mainthm400} The same statement of Theorem \ref{TH1}
remains correct when $z\in Q$ is an interior point, with $B_{x,r}^+$ replaced
by $B_{x,r}$.
\end{thm}

\section*{Acknowledgments}

{The research of the first and third authors is partly supported
by NSERC grants nos. 22R80976 and 22R81253.  The second author is
partly supported by a PIMS PDF.}

\vspace{2mm}
Stephen Gustafson\\
Department of Mathematics, the University of British Columbia\\ 
Room 121, 1984 Mathematics Road\\ 
Vancouver, B.C., Canada V6T 1Z2\\
E-mail\,\,:\,\,gustaf@math.ubc.ca\\
\\
Kyungkeun Kang\\
Department of Mathematics, the University of British Columbia\\ 
Room 121, 1984 Mathematics Road\\ 
Vancouver, B.C., Canada V6T 1Z2\\
E-mail\,\,:\,\,kkang@math.ubc.ca\\
\\
Tai-Peng Tsai\\
Department of Mathematics, the University of British Columbia\\ 
Room 121, 1984 Mathematics Road\\ 
Vancouver, B.C., Canada V6T 1Z2\\
E-mail\,\,:\,\,ttsai@math.ubc.ca\\

\end{document}